\documentclass[11pt]{article}

\newenvironment{pf}{\begin{trivlist}\item[] \textbf{Proof.}}
                     {\hspace*{\fill} $\square$\end{trivlist}}

\newenvironment{pfofpropP1}{\begin{trivlist}\item[] \textbf{Proof of Proposition \ref{P1}.}}
                     {\hspace*{\fill} $\square$\end{trivlist}}

\newenvironment{pfofpropP2}{\begin{trivlist}\item[] \textbf{Proof of Corollary \ref{P2}.}}
                     {\hspace*{\fill} $\square$\end{trivlist}}

\newenvironment{pfofthmT1}{\begin{trivlist}\item[] \textbf{Proof of Theorem \ref{T1} and Corollary \ref{C1}.}}
                     {\hspace*{\fill} $\square$\end{trivlist}}

\title{Right inverses of L\'evy processes: \\ the excursion measure in the general case} 
\author{
Mladen Savov\thanks{University of Oxford; email: savov@stats.ox.ac.uk} \and Matthias Winkel\thanks{University of Oxford; email: winkel@stats.oc.ax.uk
}
 }

\usepackage{amsfonts,amsmath,amssymb,bbm}
\usepackage{amsfonts,amsmath,latexsym,amssymb,amsthm}

\usepackage{verbatim}

\newcommand{\LL}{L\'{e}vy }

\newcommand{\beq}{\begin{eqnarray*}}
\newcommand{\eeq}{\end{eqnarray*}}
\newcommand{\ds}{\displaystyle}

\newcommand{\beqs}{\begin{eqnarray*}&\ds}
\newcommand{\eeqs}{&\end{eqnarray*}}

\newcommand{\bE}{\mathbb{E}}

\newcommand{\bK}{\mathbb{K}}

\newcommand{\bP}{\mathbb{P}}
\newcommand{\bQ}{\mathbb{Q}}
\newcommand{\bR}{\mathbb{R}}

\newcommand{\cD}{\mathcal{D}}
\newcommand{\cE}{\mathcal{E}}

\newcommand{\cL}{\mathcal{L}}

\newcommand{\cU}{\mathcal{U}}

\newtheorem{theorem}{Theorem}

\newtheorem{proposition}[theorem]{Proposition}
\newtheorem{lemma}[theorem]{Lemma}
\newtheorem{corollary}[theorem]{Corollary}

\begin{document}
\maketitle

\begin{abstract} This article is about right inverses of L\'evy processes as first introduced by Evans in 
  the symmetric case and later studied systematically by the present authors and their co-authors. Here we add to the existing fluctuation theory an
  explicit description of the excursion measure away from the (minimal) right inverse. This description unifies known
  formulas in the case of a positive Gaussian coefficient and in the bounded variation case. While these known formulas
  relate to excursions away from a point starting negative continuously, 
  and excursions started by a jump, the present description is in terms of excursions away from the supremum continued up to a
  return time. In the unbounded variation case with zero Gaussian coefficient previously excluded, excursions start negative
  continuously, but the excursion measures away from the right inverse and away from a point are mutually singular. We also 
  provide a new construction and a new formula for the Laplace exponent of the minimal right inverse.
    
\emph{AMS 2000 subject classifications: 60G51\newline
}

\end{abstract}

\section{Introduction}

Evans \cite{Ev-00} defined a \em right inverse \em of a L\'evy process $X=(X_t,t\ge 0)$ to be any increasing process $K=(K_x,x\ge 0)$ such that $X_{K_x}=x$ for all $x\ge 0$. A \em partial right inverse \em \cite{Win-02} is any increasing 
process $K=(K_x,0\le x<\xi_K)$ such that $X_{K_x}=x$ for all $0\le x<\xi_K$ for some (random) $\xi_K>0$. The existence of partial right inverses is a local path property that has been completely characterised \cite{DoSa-10,Ev-00,Win-02} in terms of the L\'evy-Khintchine triplet $(a,\sigma^2,\Pi)$ of the L\'evy process $X$, i.e. $a\in\bR$, $\sigma^2\ge 0$ and $\Pi$ measure on $\mathbb{R}$ with $\Pi(\{0\})=0$ and $\int_{\bR}(1\wedge y^2)\Pi(dy)<\infty$ such that
$$\mathbb{E}\left(e^{i\lambda X_t}\right)=e^{-t\psi(\lambda)},\qquad\mbox{where }\psi(\lambda)=-ia\lambda+\frac{1}{2}\sigma^2\lambda^2+\int_{\bR}\left(1-e^{i\lambda y}+i\lambda y1_{\{|y|\le 1\}}\right)\Pi(dy).$$ 
Where right inverses exist, the \em minimal right-continuous \em right inverse is a subordinator, and for partial right inverses (and $\xi_K$ maximal), a subordinator run up to an independent exponential time $\xi_K$. In the sequel, we focus on this minimal right-continuous (partial) right inverse and denote the Laplace exponent by
\begin{equation}\label{rho}\rho(q)=-\ln\left(\bE\left(e^{-qK_1};\xi_K>1\right)\right)=\kappa_K+\eta_K q+\int_{(0,\infty)}\left(1-e^{-qt}\right)\Lambda_K(dt).
\end{equation}
Evans \cite{Ev-00} showed 
further that the reflected process $Z=X-L$ is a strong Markov process, where $L_t=\inf\{x\ge 0:K_x>t\}$, $0\le t<K_{\xi_K}$, $L_t=\xi_K$, $t\ge K_{\xi_K}$; and $L$ is a local time process of $Z$ at zero. In analogy with the classical theory of excursions away from the supremum (see below), there is an associated excursion theory that studies the
Poisson point process $(e_x^Z,0\le x<\xi)$ of \em excursions away from the right inverse\em, where
$$e^Z_x(r)=Z_{K_{x-}+r},\quad 0\le r\le\Delta K_x=K_x-K_{x-},\qquad e^Z_x(r)=0,\quad r\ge\Delta K_x.$$
Specifically, we denote by $n^Z$ its intensity measure on the space $(E,\cE)$ of excursions 
$$E=\left\{\omega\in D:\omega(s)=0\mbox{ for all $s\ge\zeta(\omega)=\inf\{r>0:\omega(r)=0\}$}\right\}$$
equipped with the restriction sigma-algebra $\cE$ induced by the Borel sigma-algebra $\cD$ associated with Skorohod's topology on the space $D=D([0,\infty),\bR)$ of c\`adl\`ag paths $\omega\colon[0,\infty)\rightarrow\bR$. The entrance laws $n_r^Z(dy)=n^Z(\{\omega\in E:\omega(r)\in dy,\zeta(\omega)>r\})$ are characterised by their Fourier-Laplace transform
\begin{equation}\label{fourlapl}\int_0^\infty e^{-qr}\int_{\bR}e^{i\lambda y}n_r^Z(dy)dr=\frac{\rho(q)-i\lambda}{q+\psi(\lambda)}-\eta_K,
\end{equation}
see \cite{Ev-00,Win-02}. The sub-stochastic semi-group within these excursions is the usual killed semi-group 
$$P_t^\dagger(y,dz)=\bP_y^\dagger(\{\omega\in D:\omega(t)\in dz,\zeta(\omega)>t\}),\qquad\mbox{with }\bP_y^\dagger=\bP((X_{t\wedge\zeta(X)},t\ge 0)\in\,\cdot\;|X_0=y)$$
as canonical measure on $E\subset D$ of the distribution of $X$ starting from $y$ and frozen when hitting zero. 
More explicit expressions for $n^Z$ are available from \cite{Win-02} in two cases. When $\sigma^2>0$, then
\begin{equation}\label{sigpos}n^Z(d\omega)=\frac{2}{\sigma^2}n^X(d\omega;\omega(s)<0\mbox{ for all $0<s<\varepsilon$ and some $\varepsilon>0$})\end{equation}
is proportional to the intensity measure $n^X$ of excursions of $X$ away from zero restricted to those starting negative. In the bounded variation case, $\sigma^2=0$ and $\int_{\bR}(1\wedge|y|)\Pi(dy)<\infty$, we have
\begin{equation}\label{bvexc}n^Z(d\omega)=\frac{1}{b}\int_{\bR}\bP_y^\dagger(d\omega)\Pi(dy),\quad\mbox{where necessarily }b=a-\int_{\bR}y1_{\{|y|\le 1\}}\Pi(dy)>0.
\end{equation}
In the present paper, we describe $n^Z$ for a general L\'evy process that possesses a partial right
inverse. This seems to answer the final open question \cite{DoSa-10} related to the notion of right inverse of a L\'evy process. However, this
study can also be seen in the light of more general subordination \cite{Sim-99} of the form $X_{T_x}=Y_x$, where $(T,Y)$ is a bivariate L\'evy process, increasing in the $T$-component. 

To formulate our main result, we recall some classical fluctuation theory \cite{bert96,Doney-07}. With any \LL process $X$ we associate the ascending ladder time and ladder height processes $(\tau, H)$, a bivariate subordinator such that $H_x=X_{\tau_{x}}=\overline{X}_{\tau_x}$ visits all suprema $\overline{X}_t=\sup\{X_s,0\le s\le t\}$, $t\ge 0$. If $X_t\rightarrow-\infty$ as $t\rightarrow\infty$, then $\tau=(\tau_x,0\le x<\xi)$ is a subordinator run up to an exponential time $\xi$. We write the Laplace exponent of $(\tau, H)$ in L\'evy-Khintchine form as
\begin{equation}\label{LaplExp}
k(\alpha,\beta)=-\ln\left(\bE\left(e^{-\alpha \tau_{1}-\beta H_{1}};\xi>1\right)\right)=\kappa+\eta\alpha + \delta\beta+\int_{[0,\infty)^2}\left(1-e^{-\alpha s-\beta y}\right)\Lambda(ds,dy).
\end{equation}
It was shown in \cite{Win-02} that $\delta>0$ whenever there exists a partial
right inverse. \em In the sequel, we will always normalise the ascending ladder processes so that $\delta=1$. \em Also, when
partial right inverses exist, then $\bP(T_{\{z\}}<\infty)>0$ for all $z>0$, where $T_{\{z\}}=\inf\{t\ge 0:X_t=z\}$. In particular,
the $q$-resolvent measure
\[U^{q}(dz)=\int_{0}^{\infty}e^{-qt}\bP(X_t\in dz)dt\]
then admits a bounded density $u^{q}(z)$ that is continuous except possibly for a discontinuity at zero in the bounded variation case, see \cite[Theorem 43.19]{Sat-99}. Now $R=X-\overline{X}$ is a strong Markov process with $\tau$ as its inverse local time; its excursions, with added height $\Delta H_x$ at freezing,  
$$e^R_x(r)=R_{\tau_{x-}+r},\quad 0\le r<\Delta \tau_x=\tau_x-\tau_{x-},\qquad e^R_x(r)=\Delta H_x,\quad r\ge\Delta \tau_x,$$
form a Poisson point process whose intensity measure we denote by $\widetilde{n}^{R}$. For $\omega\in D$, we write $\zeta^+(\omega)=\inf\{r>0:\omega(r)>0\}$. For $\omega_1\in D$ with $\zeta^+(\omega_1)<\infty$ and $\omega_2\in D$, we 
concatenate
$$\omega=\omega_1\oplus\omega_2,\qquad\mbox{where }\omega(r)=\omega_1(r),\quad 0\le r<\zeta^+(\omega_1),\qquad \omega(\zeta^+(\omega_1)+r)=\omega_2(r),\quad r\ge 0.$$

\begin{theorem}\label{T1}
Let $X$ be a L\'evy process that possesses a partial right inverse. Then the excursion measures $n^Z$ away from the right inverse and $\widetilde{n}^R$ away from the supremum are related as 
\begin{equation}\label{genexc}n^Z(d\omega)=(\widetilde{n}^R\oplus\bK)(d\omega),
\end{equation} 
where the stochastic kernel
$$\bK(\omega_1,d\omega_2)=\bP_{\omega_1(\zeta^+(\omega_1))}^\dagger(d\omega_2),\quad\mbox{if }\zeta^+(\omega_1)<\infty,\qquad\bK(\omega_1,d\omega_2)=\delta_0\quad\mbox{otherwise},$$
associates to a path $\omega_1$ that passes positive a L\'evy process path $\omega_2$ that starts at the first positive height of $\omega_1$ and is frozen when reaching zero, and where
$(\widetilde{n}^R\oplus\bK)(d\omega)$ is the image measure of $\bK(\omega_1,d\omega_2)\widetilde{n}^R(d\omega_1)$ under concatenation $(\omega_1,\omega_2)\mapsto\omega=\omega_1\oplus\omega_2$. 
\end{theorem}

For the case $\sigma^2=0$ and $X$ of unbounded variation, it was noted in \cite{DoSa-10} that a.e. excursion under
$n^X$ starts positive and ends negative, while a.e. excursion under $\widetilde{n}^R$ starts negative; by Theorem \ref{T1}, the two measures $n^{Z}$ and $n^{X}$ are therefore singular, in contrast to (\ref{sigpos}) in the case $\sigma^2>0$. When $X$ is of bounded variation, the discussion in \cite[Section 5.3]{Win-02} easily yields the compatibility of (\ref{genexc}) and (\ref{bvexc}); it can happen that a jump in the ladder height process $\Delta H_x$ occurs without an excursion away from the supremum, $\Delta\tau_x=0$, and so $n^{Z}$ can charge paths with $\omega(0)>0$, while in the unbounded variation case we have $\omega(0)=0$ for $n^{Z}$-a.e. $\omega\in E$.

Other descriptions of $n^Z$ follow: let $\widetilde{n}_t^R(dz)=\widetilde{n}^R(\{\omega\in E:\omega(t)\in dz,\zeta^+(\omega)\wedge\zeta(\omega)>t\})$.

\begin{corollary}\label{C1} In the setting of Theorem \ref{T1}, the entrance laws of $n^Z$ are given by
\begin{equation}\label{T1E1}
n^{Z}_{t}(dz)=\widetilde{n}_t^R(dz)+\int_{[0,t]\times(0,\infty)}P_{t-s}^\dagger(y,dz)\Lambda(ds,dy),
\end{equation}
and the semi-group of $n^Z$, or rather $n^Z(\{\omega\in E:(\omega(t),0\le t<\zeta(\omega))\in\,\cdot\,\}$, is $(P_t^\dagger(y,dz),t\ge 0)$.
\end{corollary}
We also record an expression for the Laplace exponent $\rho$ of the partial right inverse $K$.

\begin{proposition}\label{P1}
Let $X$ be a \LL process which possesses a partial right inverse. Then
\begin{equation}\label{P1E1}
\rho(q)=\kappa+\eta q+\int_{[0,\infty)\times[0,\infty)}\left(1-e^{-qs}\frac{u^{q}(-y)}{u^{q}(0)}\right)\Lambda(ds,dy),\quad\mbox{with $u^q(0)=u^q(0+)$,}
\end{equation}
where $\kappa\geq0$, $\eta\geq 0$ and $\Lambda$ are as in {\rm(\ref{LaplExp})}, respectively, the killing rate and the drift coefficient of the ascending ladder time process and the \LL measure of the bivariate ladder subordinator $(\tau,H)$. In particular, the
characteristics $(\kappa_K,\eta_K,\Lambda_K)$ of $K$ in {\rm(\ref{rho})} are given by
$$\kappa_K=\kappa+\int_{[0,\infty)^2}\bP(T_{\{-y\}}=\infty)\Lambda(ds,dy),\qquad \eta_K=\eta,$$
\begin{equation}\Lambda_K(dt)=\int_{[0,\infty)\times[0,\infty)}\bP(s+T_{\{-y\}}\in dt;T_{\{-y\}}<\infty)\Lambda(ds,dy).\label{LambdaK}
\end{equation}
\end{proposition}
We stress that $\Lambda(\{0\},dy)$ is the zero measure unless $X$ can jump into its new supremum from the position of its current supremum. The latter can happen only when the ascending ladder time has a positive drift $\eta>0$, i.e. in particular only when $X$ is of bounded variation.

Let us note a simple consequence of Proposition \ref{P1} which can also be seen directly using repeated trials arguments: $\bP(\xi_K>x,K_x\le t)>0$ for some $t>0$ implies $\bP(\xi_K>x)=1$.

\begin{corollary}\label{P2}
A recurrent \LL process has a partial right inverse iff it has a full right inverse.
\end{corollary}

We proceed as follows. In Section \ref{Preliminaries} we recall Evans' construction, introduce an alternative construction of right inverses and heuristically derive Proposition \ref{P1}. Sections \ref{Proofs} and \ref{Proofsexc} contain proofs of Proposition \ref{P1} and Corollary \ref{P2}, and, respectively, of Theorem \ref{T1} and Corollary \ref{C1}.

\section{Construction of minimal partial right inverses}\label{Preliminaries}

We recall Evans' construction of right inverses. Recursively define for each $n\ge 0$ times
\[T_0^{(n)}=0,\qquad T^{(n)}_{k+1}=\inf\left\{t\geq T^{(n)}_{k}:\,X_{t}=\frac{k+1}{2^{n}}\right\},\quad k\ge 1,\]
and a process 
\[K^{(n)}_{x}=T^{(n)}_{k},\qquad\frac{k}{2^n}\leq x<\frac{k+1}{2^{n}},\qquad k\ge 0.\]
Then, a pathwise argument shows that, if the limit
\begin{equation}\label{RI}
K_x=\inf_{y>x}\sup_{n\geq 0}K^{(n)}_{y}
\end{equation}
is finite for $0\le x<\xi_K$ and $\xi_K>0$, it is the \em minimal \em partial right inverse: for all right-continuous partial right inverses $(U_x,0\le x<\xi^\prime)$, we have $\xi^\prime\le \xi_K$ and $U_x\geq K_x$ for all $0\le x<\xi^\prime$. 

Let us now suppose that $X$ possesses right inverses. We introduce an alternative construction that we use to sketch a heuristic proof of Proposition \ref{P1}. Informally, for all $n\geq 0$ we approximate $K$ by the ascending ladder time process $\tau$, but when an excursion away from the supremum with $e_x^R(\Delta\tau_x)=\Delta H_x>2^{-n}$ appears, our approximation $\widetilde{K}(n)$ of $K$ makes a jump whose size is the length of this excursion plus the time needed by $X$ to return to the starting height 
$H_{x-}=\overline{X}_{\tau_{x-}}=X_{\tau_{x-}}$ of the excursion, after which we iterate the procedure. Then $X_{\widetilde{K}_x(n)}$ evolves like an ascending ladder height process $H$ with jumps of sizes exceeding $2^{-n}$ removed.

Let us formalize this. Fix $n\geq 0$. Consider the bivariate subordinator $(\tau,H)$. 
To begin an inductive definition, let 
\[S_{1}(n)=\inf\left\{x\geq 0:\,\Delta H_{x}>2^{-n}\right\}\]
\[\widetilde{K}_{x}(n)=\tau_{x},\quad 0\le x<S_{1}(n)\]
\[\widetilde{K}_{S_1(n)}(n)=\inf\left\{t\ge\tau_{S_1(n)}:X_t=X_{\tau_{S_{1}(n)-}}\right\},\quad\mbox{if $S_1(n)<\infty$.}\]
Given $(\widetilde{K}_{x}(n),\,0\leq x\leq S_{m}(n))$ and $T_{m}(n)=\widetilde{K}_{S_{m}(n)}(n)<\infty$, let $X^{(m)}_{t}(n)=X_{T_{m}(n)+t}-X_{T_{m}(n)}$, $t\ge 0$. With $(\tau^{(m)}(n),H^{(m)}(n))$ as the bivariate ladder subordinator of $X^{(m)}(n)$, define 
\[S_{m+1}(n)=S_m(n)+\inf\left\{x\geq 0:\,\Delta H^{(m)}_{x}(n)>2^{-n}\right\}\]
\[\widetilde{K}_{S_{m}(n)+x}(n)=T_{m}(n)+\tau^{(m)}_{x}(n),\quad 0\leq x< S_{m+1}(n)-S_m(n)\]
\[\widetilde{K}_{S_{m+1}(n)}(n)=T_{m}(n)+\inf\left\{t\geq\tau^{(m)}_{S_{m+1}(n)-S_m(n)}(n):X^{(m)}_t(n)=X^{(m)}_{\tau^{(m)}_{S_{m+1}(n)-S_m(n)}(n)-}(n)\right\}.\]
Thus we have defined $\widetilde{K}_{x}(n)$ for all $x\ge 0$ and $n\ge 0$ a.s. Now it \em must be expected that \em
\begin{equation}K_x=\inf_{y>x}\sup_{n\ge 0}\widetilde{K}_{L_n(y)}(n)=\lim_{n\rightarrow\infty}\widetilde{K}_x(n),\quad\mbox{where }L_n(y)=\inf\{x\ge 0:X_{\widetilde{K}_x(n)}\ge y-2^{-n}\}.\label{toshow}
\end{equation}
%
To derive formula (\ref{LambdaK}) of Proposition \ref{P1}, consider the Poisson point process $(\Delta\widetilde{K}_x(n),x\ge 0)$ whose intensity measure we can calculate from $((\Delta\tau_x,\Delta H_x),x\ge 0)$ using standard thinning (keep if $\Delta H_x\le 2^{-n}$, modify if $\Delta H_x>2^{-n}$), marking by independent $T_{\{-\Delta H_x\}}$ if $\Delta H_x>2^{-n}$ and mapping $(\Delta\tau_x,\Delta H_x,T_{\{-\Delta H_x\}})\mapsto\Delta\tau_x+T_{\{-\Delta H_x\}}=\Delta\widetilde{K}_x(n)$ of Poisson point processes, as
$$\Lambda_{\widetilde{K}(n)}(dt)=\int_{y\in[0,2^{-n}]}\Lambda(dt,dy)+\int_{(s,y)\in[0,\infty)\times(2^{-n},\infty)}\bP\left(s+T_{\{-y\}}\in dt\right)\Lambda(ds,dy),
$$
which converges to the claimed expression, as $n\rightarrow\infty$. We make this approach rigorous in Appendix \ref{appx}, the main task being to rigorously establish a variant of (\ref{toshow}). In the next section, we will instead start from Evans' construction (\ref{RI}) and exploit recent developments \cite{DoK-06} on joint laws of first passage variables, which follows on more naturally from 
previous work.


\section{The Laplace exponent of $K$; proof of Proposition \ref{P1}}\label{Proofs}

Before we formulate and prove some auxiliary results, we introduce some notation. Denote by 
$\cU(ds,dy)=\int_0^\infty\bP(\tau_x\in ds,H_x\in dy;\xi>x)dx$ the potential measure of the bivariate ascending ladder subordinator
$(\tau,H)$, by $T_x^+=\inf\{t\ge 0:X_t\in(x,\infty)\}$ the first passage time across level $x>0$, by 
$\overline{X}_t=\sup_{0\le s\le t}X_s$ the supremum process, by $\overline{G}_t=\sup\{s\le t:X_s=\overline{X}_t\mbox{ or }X_{s-}=\overline{X}_t\}$ the time of the last visit to the supremum and by $O_x=X_{T_x^+}-x\ge 0$ the overshoot over level $x$. Then on $\{(t,s,w,y):t\ge 0,s\ge 0,w>0,0\le y\le x\}$, we have
  \begin{equation}\label{quintuple}
  \bP(T_x^+-\overline{G}_{T_x^+-}\in dt,\overline{G}_{T_x^+-}\in ds,O_x\in dw,x-\overline{X}_{T_x^+-}\in dy;T_x^+<\infty)=\Lambda(dt,du+y)\cU(x-dy,ds),
  \end{equation}
by a corollary of the quintuple law of Doney and Kyprianou \cite{DoK-06}.
   
Suppose that $X$ possesses a partial right inverse. Recall construction \eqref{RI}. A crucial quantity there is the hitting time of levels, $T_{\{x\}}$. A key observation for our developments is that 
\begin{equation}\label{tsplit}T_{\{x\}}=T_x^{+}+\widetilde{T}_{\{-O_x\}}\qquad\mbox{a.s. on $\{T_x^+<\infty\}$},\end{equation}
where $\widetilde{T}_{\{-O_x\}}=\inf\{t\geq 0:\widetilde{X}_{t}=-O_x\}$ for $\widetilde{X}=(X_{T_x^++t}-X_{T_x^+},t\ge 0)$ independent of $(T_x^+,O_x)$. For $q>0$, let
$$z_n=2^n\bE\left(1-e^{-qT_{\{2^{-n}\}}}\right)=2^n\bE\left(1-e^{-q(T_{2^{-n}}^++\widetilde{T}_{\{-O_{2^{-n}}\}})}\right),$$
with the convention that $e^{-\infty}=0$ and $T_{2^{-n}}^++\widetilde{T}_{\{O_{2^{-n}}\}}=\infty$ on $\{T_{2^{-n}}^+=\infty\}$.
As was already exploited by Evans \cite{Ev-00}, construction \eqref{RI} allows us to express  
\begin{equation}
\label{E20}
\rho(q)=-\ln\left(\lim_{n\rightarrow\infty}\bE\left(e^{-qK_{1}^{(n)}}\right)\right)=-\ln\left(\lim_{n\rightarrow\infty}\left(1-\frac{z_n}{2^n}\right)^{2^n}\right)=\lim_{n\rightarrow\infty}z_n,
\end{equation}
because $K_1^{(n)}$ is the sum of $2^n$ independent random variables with the same distribution as $T_{\{2^{-n}\}}$. To calculate this limit, we will use (\ref{tsplit}) and also decompose $z_n$, as follows, setting
\begin{equation}\label{E21}\widehat{z}_n=2^n\bE\left(\left(1-e^{-qT_{\{2^{-n}\}}}\right)1_{\{O_{2^{-n}}>0,T_{2^{-n}}^+<\infty\}}\right)\quad\mbox{and}\quad\widetilde{z}_n=z_n-\widehat{z}_n.
\end{equation}

\begin{lemma}\label{L1}
Let $X$ be a \LL process of unbounded variation which possesses a partial right inverse. Then
\begin{equation}\label{L1E1}
\lim_{n\to\infty}\widehat{z}_{n}=\int_{(t,h)\in(0,\infty)^2}\left(1-e^{-qt}\frac{u^q(-h)}{u^q(0)}\right)\Lambda(dt,dh).
\end{equation}
\end{lemma} 
\begin{pf} 
According to \cite[Theorems 43.3, 43.19 and 47.1]{Sat-99}, the resolvent density $u^{q}$ is bounded and continuous for all $q>0$ and
$\bE (e^{-qT_{\{x\}}})=u^{q}(x)/u^{q}(0)$ for all $x\in\bR$. We use (\ref{quintuple}) to obtain
\beq&&\hspace{-0.5cm}\bE\left(\left(1-e^{-q(T_x^{+}+\widetilde{T}_{\{-O_{x}\}})}\right)1_{\{O_{x}>0\}}\right)=\bE\left(\left(1-e^{-q(\overline{G}_{T_x^{+}-}+T_x^{+}-\overline{G}_{T_x^{+}-}+\tilde{T}_{\{-O_x\}})}\right)1_{\{O_x>0\}}\right)\\
       &&=\int_{(s,y)\in(0,\infty)\times[0,x]}\int_{(t,w)\in(0,\infty)^2}\left(1-e^{-q(s+t)}\frac{u^q(-w)}{u^q(0)}\right)\Lambda(dt,dw+y)\cU(ds,x-dy)\\
       &&=\int_{(t,h)\in(0,\infty)^2}\int_{(s,y)\in(0,\infty)\times[0,x\wedge h]}\left(1-e^{-q(s+t)}\frac{u^q(-h+y)}{u^q(0)}\right)\cU(ds,x-dy)\Lambda(dt,dh).
  \eeq
Therefore it will be sufficient to show that
as $n$ tends to infinity, we have the convergence
    \begin{eqnarray}\label{E1} \nonumber&&\hspace{-0.5cm}2^n\int_{(t,h)\in(0,\infty)^2}\int_{(s,y)\in(0,\infty)\times[0,2^{-n}\wedge h]}\left(1-e^{-q(t+s)}\frac{u^q(-h+y)}{u^q(0)}\right)\cU(ds,2^{-n}-dy)\Lambda(dt,dh)\\
		&&\longrightarrow \int_{(t,h)\in(0,\infty)^2}\left(1-e^{-qt}\frac{u^q(-h)}{u^q(0)}\right)\Lambda(dt,dh).
	\end{eqnarray}
First fix $(t,h)\in(0,\infty)^2$ and consider the bounded and continuous function $f:[0,\infty)^2\rightarrow[0,\infty)$
    given by
    \begin{equation}\label{func}
    f(s,y)=1-e^{-q(t+s)}\frac{u^q(-h+y)}{u^q(0)}
    \end{equation}
    and the measures $\vartheta_{n,h}(ds,dy)=1_{[0,2^{-n}\wedge h]}(y)2^n\cU(ds,2^{-n}-dy)$ on $[0,\infty)^2$. Since $H$ has unit drift, we have that $\bP(H_t\ge t\mbox{ for all $t\ge 0$})=1$ and so for all $\varepsilon>0$, we have
    \begin{align*}   &\lim_{n\to\infty}\vartheta_{n,h}\left(([0,\varepsilon]\times[0,\varepsilon])^{c}\right)=\lim_{n\to\infty}2^n\int_0^{\infty}\bP\left(\tau_x>\varepsilon, H_x\leq 2^{-n}\right)dx\\
    &=\lim_{n\to\infty}2^n\int_0^{2^{-n}}\bP\left(\tau_x>\varepsilon, H_x\leq 2^{-n}\right)dx\leq\lim_{n\to\infty}\bP(\tau_{2^{-n}}>\varepsilon)=0,
    \end{align*}
    whereas $H_t/t\rightarrow 1$ a.s., as $t\rightarrow 0$, implies that $\bP(H_{2^{-n}(1-\varepsilon)}\le 2^{-n})\ge 1-\varepsilon$ for $n$ sufficiently large, and so
    $$1\ge 2^n\int_0^{2^{-n}}\bP(\tau_x\ge 0,H_x\le 2^{-n})dx\ge 2^n\int_0^{2^{-n}(1-\varepsilon)}\bP(H_x\le 2^{-n})dx\ge (1-\varepsilon)^2.$$
    This shows convergence $1_{[0,2^{-n}\wedge h]}(y)2^n\cU(2^{-n}-dy,ds)\rightarrow\delta_{(0,0)}$ weakly as $n$ tends to 
    infinity, where $\delta_{(0,0)}$ is the Dirac measure in $(0,0)$. This shows that 
    \begin{equation}\label{E7}
    \int_{[0,\infty)^{2}}f(s,y)\vartheta_{n,h}(ds,dy)\rightarrow f(0,0).
    \end{equation}
    
    To deduce \eqref{E1}, and hence \eqref{L1E1}, from \eqref{E7}, we use the Dominated Convergence Theorem and for this purpose we show that
    \begin{equation}\label{fn}f_{n}(t,h)=\int_{(s,y)\in(0,\infty)\times[0,2^{-n}\wedge h]}\left(1-e^{-q(t+s)}\frac{u^q(-h+y)}{u^q(0)}\right)\cU(ds,2^{-n}-dy)\end{equation}
    is bounded above by a $\Lambda$-integrable function. We shall prove the bound
     \begin{equation}\label{E2} f_{n}(t,h)
      \le(1-e^{-qt})+\left(k(q,0)+\frac{1}{2}k(0,\rho(q))\right)(1\wedge h)+\left(1-\frac{u^q(-h)}{u^q(0)}\right),
  \end{equation}
  where we recall from (\ref{LaplExp}) that $k$ is the Laplace exponent of $(\tau,H)$. The integrand in \eqref{fn} is
    \beq f(s,y)&=&(1-e^{-qt})+e^{-qt}(1-e^{-qs})+e^{-q(t+s)}\left(1-\frac{u^q(-h+y)}{u^q(0)}\right)\\
			&\le&(1-e^{-qt})+(1-e^{-qs})+\left(1-\frac{u^q(-h+y)}{u^q(0)}\right),
	\eeq
	three terms, where $f(s,y)$ is defined in \eqref{func}.
	First note that as before since $H$ has unit drift
	\[\vartheta_{n,h}([0,\infty)^{2})=2^{n}\int_{0}^{2^{-n}}P(H_{t}\leq 2^{-n})dt\leq 1.\]
	Therefore
	\begin{equation}\label{E3}
	\int_{(s,y)\in [0,\infty)^{2}}(1-e^{-qt})\vartheta_{n,h}(ds,dy)\leq 1-e^{-qt}.
	\end{equation}
For the second term we use $\bP(H_u-H_t\ge u-t\mbox{ for all $u\ge t\ge 0$})=1$ to write
	\begin{eqnarray}\label{E4} \nonumber&&\hspace{-0.5cm}\int_{[0,\infty)^{2}}(1-e^{-qs})\vartheta_{n,h}(ds,dy)\\
	     \nonumber&&\le\int_{s\in[0,\infty)}(1-e^{-qs})2^n\int_{r\in[0,\infty)}\bP\left(\tau_{r}\in ds,H_r\in[2^{-n}-h,2^{-n}]\right)dr\\
	     \nonumber&&=2^n\int_{[0,2^{-n}]}\bE\left(\left(1-e^{-q\tau_{r}}\right)1_{\{H_r\in[2^{-n}-h,2^{-n}]\}}\right)dr\\
	    \nonumber &&\le 2^n\int_{[0,2^{-n}]}\bE\left(\left(1-e^{-q\tau_{2^{-n}}}\right)1_{\{H_r\in[2^{-n}-h,2^{-n}]\}}\right)dr\\
	     &&\le 2^n(h\wedge 1)\left(1-e^{-2^{-n}k(q,0)}\right)\le(h\wedge 1)k(q,0).
	\end{eqnarray}
	 For the third term we mimick the previous calculation to get
    \beq &&\hspace{-0.5cm}\int_{[0,\infty)^{2}}\left(1-\frac{u^q(-h+y)}{u^q(0)}\right)\vartheta_{n,h}(ds,dy)\\
		 &&=\bE\left(2^n\int_{r\in[0,2^{-n}\wedge h]}\left(1-\frac{u^q((-h+H_r)\wedge 0)}{u^q(0)}\right)dr\right)=:\Upsilon(h).
    \eeq
    We now exploit the fact \cite[Corollary 2]{Win-02} that $x\mapsto e^{-\rho(q)x}u^q(-x)$ is decreasing, and also 
    $1-e^{-x}\ge x-x^2/2$, to see, for $h\le 2^{-n}$
    \beq \Upsilon(h)&\le&2^n\left(h-\frac{u^q(-h)}{u^q(0)}\bE\left(\int_0^h e^{-\rho(q)H_r}dr\right)\right)\\
					  &=&2^n\left(h-\frac{u^q(-h)}{u^q(0)}\frac{1}{k(0,\rho(q))}(1-e^{-h k(0,\rho(q))})\right)\\
					  &\le&\left(1-\frac{u^q(-h)}{u^q(0)}\right)+h\frac{1}{2}k(0,\rho(q)),
    \eeq
    and similarly for $h>2^{-n}$,
    $$ \Upsilon(h)\le 2^n\left(2^{-n}-\frac{u^q(-h)}{u^q(0)}\left(2^{-n}+2^{-2n}\frac{k(0,\rho(q))}{2}\right)\right).
    $$
    Together, this yields an upper bound for all $h\in(0,\infty)$
    \begin{equation}\label{E5} \Upsilon(h)\le\left(1-\frac{u^q(-h)}{u^q(0)}\right)+\frac{1}{2}k(0,\rho(q))(h\wedge 1).\end{equation}
    
    Thus \eqref{E2} follows from \eqref{E3}, \eqref{E4} and \eqref{E5}. In view of the fact that $\Lambda(dt,dh)$ is a \LL measure of a subordinator the RHS of \eqref{E2} will be $\Lambda(dt,dh)$-integrable if $1-u^{q}(-h)/u^{q}(0)$ is $\Lambda(dt,dh)$-integrable. First using Fubini's Theorem in \eqref{E1}, followed by Fatou's Lemma, because of \eqref{E7} and the simple inequality
    \[1-\frac{u^{q}(-h)}{u^{q}(0)}\leq 1-e^{-qt}\frac{u^{q}(-h)}{u^{q}(0)},\]
    we get 
    \[\liminf_{n\to\infty}\hat{z}_{n}\geq \int_{(t,h)\in(0,\infty)^2 }\left(1-\frac{u^{q}(-h)}{u^{q}(0)}\right)\Lambda(dt,dh).\]
    On the other hand, \eqref{E20} gives $\lim_{n\to\infty}z_{n}=\rho(q)<\infty$. Moreover, $\widehat{z}_{n}\leq z_{n}$ and we conclude that $1-u^{q}(-h)/u^{q}(0)$ is $\Lambda(dt,dh)$-integrable. Thus \eqref{E2}, together with the Dominated Convergence Theorem, implies \eqref{E1} and then \eqref{L1E1}.
\end{pf}

 \begin{lemma}\label{L2}
 Assume that $X$ is of unbounded variation and that $X$ possesses a partial right inverse. Then
$$ \lim_{n\to\infty}\widetilde{z}_{n}=\kappa+\int_{(0,\infty)}(1-e^{-qu})\Lambda(dt,\{0\}).$$
 \end{lemma}   
 \begin{pf} This proof is based on \cite[Theorem VI.18]{bert96}, which yields  $\mathbb{E}(1-e^{-q T_x^+})=k(q,0)V^q(x)$, where
  $V^q(x)=\int_0^\infty\mathbb{E}(e^{-q\tau_s};H_s\le x)ds$. Also, $V^q(x)\sim x$ as $x\downarrow 0$, since
  $V^q$ is differentiable with $v^q(0+)=1$, by dominated convergence, as $H$ has unit drift coefficient. Then note that 
  \beq\widetilde{z}_n&=&2^n\bE\left(\left(1-e^{-qT_{\{2^{-n}\}}}\right)1_{\{O_{2^{-n}}=0,T_{2^{-n}}^+<\infty\;{\rm or}\;T_{2^{-n}}^+=\infty\}}\right)\\ &=&2^n\bE\left(\left(1-e^{-qT_{2^{-n}}^+}\right)1_{\{O_{2^{-n}}=0,T_{2^{-n}}^+<\infty\;{\rm or}\;T_{2^{-n}}^+=\infty\}}\right)\eeq
  and so, we obtain the required formula from \eqref{LaplExp} noting $\eta=0$ in the unbounded variation case
  \begin{eqnarray*}&&\hspace{-0.5cm}x^{-1}\mathbb{E}\left(1-e^{-q T_x^+}\right)-x^{-1}\mathbb{E}\left(\left(1-e^{-qT_x^+}\right)1_{\{O_{x}>0,T_x^+<\infty\}}\right)\\
	&&\longrightarrow \kappa+\int_{(t,h)\in(0,\infty)\times[0,\infty)}(1-e^{-qt})\Lambda(dt,dh)-\int_{(t,h)\in(0,\infty)^2}(1-e^{-q t})\Lambda(dt,dh).
  \end{eqnarray*}
\end{pf}
Although this is not necessary for our proof of Lemma \ref{L2}, we would like to mention that we can calculate explicitly  $\widetilde{z}_n$ or, as Andreas Kyprianou pointed out to us,   
$$\bP(T_x^+\in dt,O_x=0)dx=\bP(T_x^+=\overline{G}_{T_x^+-}\in dt,O_x=0,x-\overline{X}_{T_x^+-}=0;T_x^+<\infty)dx=\cU(dt,dx),$$
which complements (\ref{quintuple}). 
Indeed, note that for $H^{-1}_x=\inf\{s\ge 0:H_s>x\}$ we have $O_x=0$ iff $\Delta H_{H_x^{-1}}=0$ and $T_x^+=\tau_{H_x^{-1}}$ for a.e. $x\ge 0$ a.s., so that as $H$ has unit drift coefficient,
\beq &&\hspace{-0.5cm}\int_{x\in(0,\infty)}\int_{t\in(0,\infty)}e^{-\alpha t-\beta x}\bP(T_x^+\in dt,O_x=0)dx=\bE\left(\int_0^\infty e^{-\alpha\tau_{H_x^{-1}}-\beta H_{H_x^{-1}}}1_{\{\Delta H_{H^{-1}_x}=0\}}dx\right)\\
&&=\bE\left(\int_0^\infty e^{-\alpha\tau_{H_x^{-1}}-\beta H_{H_x^{-1}}}dH_x^{-1}\right)=\bE\left(\int_0^\infty e^{-\alpha\tau_{s}-\beta H_{s}}ds\right)=\int_{(x,t)\in(0,\infty)^2}e^{-\alpha t-\beta x}\cU(dt,dx).
\eeq
Next we prove Proposition \ref{P1}.

\begin{pfofpropP1} Let $X$ be a L\'evy process that possesses a partial right inverse. If $X$ is of unbounded variation, we 
have $\lim_{n\to \infty}z_{n}=\rho(q)$ by \eqref{E20}; and \eqref{E21} together with Lemma \ref{L1} and Lemma \ref{L2} proves the claim. 

When $X$ is of bounded variation we refer to \cite[Section 5.3]{Win-02} which discusses in this case the Laplace exponent 
$\rho(q)$ of the partial right inverse and how the right inverse relates to ladder processes.

This establishes (\ref{P1E1}). The characteristics can now be read off by inverting the Laplace transform $u^q(-y)/u^q(0)=\bE(e^{-qT_{\{-y\}}};T_{\{-y\}}<\infty)$ for $y\ge 0$, where we recall $u^q(0)=u^q(0+)$, which entails
$\bP(T_{\{0\}}=0)=1$. 
\end{pfofpropP1}

Let us briefly explore the context of the last part of this proof. For $y=0$, note that $\bP(T_{\{0\}}=0)=1$ is a trivial consequence of the definition $T_{\{0\}}=\inf\{t\ge 0:X_t=0\}$. Indeed, this is the appropriate notion to use in the light of Theorem \ref{T1}, where excursions below the supremum ending at zero before passing positive do \em not \em get marked by a
further return time $T^{>}_{\{0\}}=\inf\{t>0:X_t=0\}$, which in the bounded variation case would have Laplace transform 
$\bE(e^{-qT_{\{0\}}^>};T_{\{0\}}^><\infty)=u^q(0-)/u^q(0+)<1$, see e.g. \cite[Theorem 43.21]{Sat-99}.  

\begin{pfofpropP2} Suppose that $X$ has a partial right inverse. In terms of the characteristics $(\kappa_K,\eta_K,\Lambda_K)$ of the minimal partial right inverse $K$, this is indeed a full right inverse if $\rho(0+)=\kappa_K=0$. But if $X$ is 
recurrent (and has a partial right inverse), then $X$ does not drift to $-\infty$, so $\kappa=0$, and $\bP(T_{\{-y\}}=\infty)=0$   for all $y\in\mathbb{R}$, so indeed 
$$\kappa_K=\kappa+\int_{[0,\infty)^2}\bP(T_{\{-y\}}=\infty)\Lambda(ds,dy)=0.$$ 
\end{pfofpropP2}

Similarly, it is also straightforward to show that in the transient case $X$ possesses a full right inverse if and only if $X$ drifts to $+\infty$ and has no positive jumps in that $\Pi((0,\infty))=0$.

\section{The excursion measure away from $K$; proof of Theorem \ref{T1}}\label{Proofsexc}

Although Theorem \ref{T1} is more refined than Proposition \ref{P1}, it is now a straightforward consequence.

\begin{pfofthmT1} The proof relies crucially on \cite[Theorem 2]{Win-02}. First recall that the semigroup within excursions away from the right inverse is $(P^\dagger_t(y,dz),t\ge 0)$. With this, the excursion measure $n^Z$ is uniquely determined by its entrance laws. Let us first check that the entrance laws given in (\ref{T1E1}) satisfy (\ref{fourlapl}), which is (vi) of \cite[Theorem 2]{Win-02}.
A standard excursion measure computation similar to the proof of that theorem together with the Wiener-Hopf factorization gives directly
\begin{eqnarray}\label{E11}
\nonumber&&\hspace{-0.5cm}\int_{0}^{\infty}e^{-qt}\int_{\mathbb{R}}e^{i\lambda y}\widetilde{n}_t^R(dy)dt=\frac{k(q,0)}{q}\int_{\mathbb{R}\setminus\{0\}}e^{i\lambda x}\bP(R_{\gamma(q)}\in dx)\\
&&=\frac{k(q,0)}{q}\left(\bE(e^{i\lambda R_{\gamma(q)}})-\bP(R_{\gamma(q)}=0)\right)
=\frac{k(q,0)\widehat{k}(q,0)}{q\widehat{k}(q,i\lambda)}-\eta=
\frac{k(q,-i\lambda)}{q+\psi(\lambda)}-\eta,
\end{eqnarray}
where $\widehat{k}(\alpha,\beta)$ is the Laplace exponent of the bivariate descending ladder process, which is the ascending
ladder process of $-X$, and where $\gamma(q)$ is an independent exponential random variable with rate parameter $q$; we recall $R_{\gamma(q)}\stackrel{d}=\underline{X}_{\gamma(q)}=\inf\{X_{s},0\le s\le\gamma(q)\}$; the Wiener-Hopf identity
\[\frac{k(q,0)}{k(q,-i\lambda)}\;\frac{\widehat{k}(q,0)}{\widehat{k}(q,i\lambda)}=\frac{q}{q+\psi(\lambda)}\]
can be found in \cite[Formulas (4.3.4) and (4.3.7)]{Doney-07}.

Next we compute the joint transform of the remaining part of the RHS of \eqref{T1E1}. 
\begin{eqnarray*}
&&\hspace{-0.5cm}\int_{0}^{\infty}e^{-qt}\int_{[0,t]\times[0,\infty)}\int_{\mathbb{R}}e^{i\lambda z}P^\dagger_{t-s}(y,dz)\Lambda(ds,dy)dt\\
&&=\int_{[0,\infty)\times(0,\infty)}e^{-qs}\int_s^\infty e^{-q(t-s)}\mathbb{E}_y(e^{i\lambda X_{t-s}};T_{\{0\}}>t-s)dt\Lambda(ds,dy)\\
&&=\frac{1}{q}\int_{[0,\infty)\times(0,\infty)}e^{-qs}\mathbb{E}_y(e^{i\lambda X_{\gamma(q)}};T_{\{0\}}>\gamma(q))\Lambda(ds,dy).
\end{eqnarray*}
where $\gamma(q)$ is an independent exponential random variable with rate parameter $q$. Also recall
\[\bP_{y}(T_{\{0\}}\leq \gamma(q))=\bE(e^{-qT_{\{-y\}}})=\frac{u^{q}(-y)}{u^{q}(0)},\qquad\mbox{where $u^q(0)=u^q(0+)$,}\]
so that
\begin{eqnarray*}
&&\hspace{-0.5cm}\frac{qe^{i\lambda y}}{q+\psi(\lambda)}=\bE_{y}(e^{i\lambda X_{\gamma(q)}})=\bE_{y}(e^{i\lambda X_{\gamma(q)}};T_{\{0\}}>\gamma(q))+\bE_{y}(e^{i\lambda X_{\gamma(q)}};T_{\{0\}}\leq \gamma(q))\\
&&=\bE_{y}(e^{i\lambda X_{\gamma(q)}};T_{\{0\}}>\gamma(q))+\bP_{y}(T_{\{0\}}\leq\gamma(q))\bE(e^{i\lambda X_{\gamma(q)}})\\
&&=
\bE_{y}(e^{i\lambda X_{\gamma(q)}};T_{\{0\}}>\gamma(q))+\frac{u^{q}(-y)}{u^{q}(0)}\frac{q}{q+\psi(\lambda)}.
\end{eqnarray*}
Thus we have computed $\bE_{y}(e^{i\lambda X_{\gamma(q)}};T_{\{0\}}>\gamma(q))$. Hence
\begin{eqnarray*}
&&\hspace{-0.5cm}\int_{0}^{\infty}e^{-qt}\int_{[0,t]\times[0,\infty)}\int_{\mathbb{R}}e^{i\lambda z}P^\dagger_{t-s}(y,dz)\Lambda(ds,dy)dt\\
&&=\frac{1}{q+\psi(\lambda)}\int_{[0,\infty)\times(0,\infty)}e^{-qs}\left(e^{i\lambda y}-\frac{u^q(-y)}{u^q(0)}\right)\Lambda(ds,dy).
\end{eqnarray*}
Next observe that by using \eqref{LaplExp} and adding the numerator of the first term in \eqref{E11} we get
\begin{eqnarray*}
&&\hspace{-0.5cm}k(q,-i\lambda)+\int_{[0,\infty)\times(0,\infty)}e^{-qs}\left(e^{i\lambda y}-\frac{u^q(-y)}{u^q(0)}\right)\Lambda(ds,dy)\\
&&=\kappa+\eta q-i\lambda+\int_{[0,\infty)^2}\left(1-e^{-qs+i\lambda y}+e^{-qs+i\lambda y}-e^{-qs}\frac{u^q(-y)}{u^q(0)}\right)\Lambda(ds,dy)=\rho(q)-i\lambda.
\end{eqnarray*}
Together with \eqref{E11}, this proves \eqref{T1E1} since \eqref{fourlapl} holds for the measure on the RHS of \eqref{T1E1}.

To finish the proof of Theorem \ref{T1}, we note that the RHS of (\ref{genexc}) is Markovian with semi-group $(P_t^\dagger(y,dz),t\ge 0)$ and check that the RHS of (\ref{genexc}) has as entrance laws the RHS of (\ref{T1E1}).
Specifically, 
\begin{eqnarray*}&&\hspace{-0.5cm}(\widetilde{n}^R\oplus\bK)(\{\omega\in D:\omega(t)\in dz;\zeta(\omega)>t\})\\
&&=\widetilde{n}^R(\{\omega_1\in D:\omega_1(t)\in dz;\zeta^+(\omega_1)>t\})\\
&&\hspace{1cm}+\int_{\{\omega_1\in D:\zeta^+(\omega_1)\le t\}}\bK(\omega_1;\{\omega_2\in D:\omega_2(t-\zeta^+(\omega_1))\in dz;\zeta(\omega_2)>t-\zeta^+(\omega_1)\})\widetilde{n}^R(d\omega_1)\\
&&=\widetilde{n}_t^R(dz)+\int_{[0,t]\times[0,\infty)}P_{t-s}^\dagger(y,dz)\Lambda(ds,dy)
\end{eqnarray*}
since $\widetilde{n}^R(\{\omega_1\in D:\zeta^+(\omega_1)\in ds,\omega_1(\zeta^+(\omega_1))\in dy)=\Lambda(ds,dy)$.
\end{pfofthmT1}

\begin{appendix}
\section{Alternative proof of Proposition \ref{P1}}\label{appx}

To simplify notation, let us assume in the sequel that $X$ possesses a full right inverse. In the case where only partial right inverses exist we can follow the second construction of Section \ref{Preliminaries} until the first $m\ge 0$ for which $H^{(m)}(n)$ is killed before its first
jump of size exceeding $2^{-n}$. If we denote the resulting process by $(\widetilde{K}_x(n),0\le x<\widetilde{\xi}(n))$, we can insert suitable restrictions to events such as $\{\widetilde{\xi}(n)>x\}$ into the following arguments.

\begin{lemma}\label{lm7} Let $\widetilde{K}_x(n)$ be as in Section {\rm \ref{Preliminaries}}. Then 
  $((\widetilde{K}_x(n),X_{\widetilde{K}_x(n)}),x\ge 0)$ is a bivariate subordinator with drift coefficient $(\eta,1)$ and 
  L\'evy measure
  $$\widetilde{\Lambda}_n(dt,dz)=\Lambda(dt,dz\cap[0,2^{-n}])+\int_{(s,y)\in[0,\infty)\times(2^{-n},\infty)}\bP(s+T_{\{-y\}}\in dt,0\in dz)\Lambda(ds,dy).$$
\end{lemma}
\begin{pf} 
  Let $\widetilde{H}_x(n)=X_{\widetilde{K}_x(n)}$. Then $(\widetilde{K}(n),\widetilde{H}(n))$ inherits the drift 
  coefficient $(\eta,1)$ from $(\tau,H)$. By standard thinning properties of Poisson point processes,   
  $((\Delta\widetilde{K}_x(n),\Delta\widetilde{H}_x(n)),0\le x<S_1(n))$ has the
  distribution of a Poisson point process with intensity measure $\Lambda(dt,dz\cap[0,2^{-n}])$ run up to an independent 
  exponential time $S_1(n)$ with parameter $\lambda=\Lambda([0,\infty)\times(2^{-n},\infty))$, and 
  $$\bP(\Delta\widetilde{K}_{S_1(n)}\in dt,\Delta\widetilde{H}_{S_1(n)}\in dz)=\lambda^{-1}\int_{[0,\infty)\times(2^{-n},\infty)}\bP(s+T_{\{-y\}}\in dt,0\in dz)\Lambda(ds,dy).$$
  By the strong Markov property of $X$ at $T_m(n)$, $m\ge 1$, the process $((\Delta\widetilde{K}_x(n),\Delta\widetilde{H}_x(n)),x\ge 0)$ 
  with points at $S_m(n)$, $m\ge 1$, removed, is a Poisson point process with intensity measure $\Lambda(dt,dz\cap[0,2^{-n}])$
  independent of the removed points, which we collect in independent and identically distributed vectors
  $(S_m(n)-S_{m-1}(n),\Delta\widetilde{K}_{S_m(n)}(n),\Delta\widetilde{H}_{S_m(n)}(n))$, $m\ge 1$. By standard superposition of
  Poisson point processes, the result follows.
\end{pf}

\begin{lemma} With $\widetilde{K}_x(n)$ as in Section {\rm \ref{Preliminaries}}, we have
\[K_x=\inf_{y>x}\sup_{n\ge 0}\widetilde{K}_{L_n(y)}(n),\qquad\mbox{where }L_n(y)=\inf\{x\ge 0:X_{\widetilde{K}_x(n)}\ge y-2^{-n}\}.\]
\end{lemma}
\begin{pf} By construction, the process $(X_{\widetilde{K}_x(n)},x\ge 0)$ has no jumps of size exceeding $2^{-n}$, so that
  $x-2^{-n}\le X_{\widetilde{K}_{L_n(x)}(n)}\le x$. Note that we have $\widetilde{K}_{L_n(x)}(n)\le K_x$. Let us define $\widetilde{K}_x$ by
  $$\widetilde{K}_x=\liminf_{n\rightarrow\infty}\widetilde{K}_{L_n(x)}(n)=\lim_{n\rightarrow\infty}\inf_{m\ge n}\widetilde{K}_{L_m(x)}(m)\le\sup_{n\ge 0}\widetilde{K}_{L_n(x)}(n)\le K_x,$$
  where the limit in the middle member of this sequence of inequalities is an increasing limit of stopping times. Since $X$ is right-continuous and quasi-left-continuous \cite[Proposition I.7]{bert96}, we obtain
  $$x-2^{-n}\le X_{\inf_{m\ge n}\widetilde{K}_{L_m(x)}}\le x\qquad\Rightarrow\qquad X_{\widetilde{K}_x}=\lim_{n\rightarrow\infty}X_{\inf_{m\ge n}\widetilde{K}_{L_m(x)}}=x\qquad\mbox{a.s.}$$
  Now, it is standard to argue that $X_{\widetilde{K}_q}=q$ holds a.s. simultaneously for all $q\in\bQ\cap[0,\infty)$
  and, since $x\mapsto\widetilde{K}_x$ is increasing and $X$ right-continuous, 
  $\inf_{y>x}\widetilde{K}_y=\inf_{q\in\bQ\cap(x,\infty)}\widetilde{K}_q\le K_x$ is a right-continuous right inverse.
  Since $K$ is the minimal right-continuous right inverse, $\inf_{y>x}\widetilde{K}_y=K_x$.  
\end{pf}

\begin{lemma}\label{lm9} Let $\widetilde{K}_x(n)$ be as in Section {\rm \ref{Preliminaries}}. Then for all $x\ge 0$ there is convergence
  along a subsequence $(n_k)_{k\ge 0}$ of $\lim_{k\rightarrow\infty}\widetilde{K}_x(n_k)=K_x$ a.s.
\end{lemma}
\begin{pf} Denote by $\cL eb$ Lebesgue measure on $[0,\infty)$. By Lemma \ref{lm7}, $\widetilde{H}(n)$ has drift coefficient
  $1$, so $\widetilde{H}_x(n)\ge x$ and $\cL eb(\{\widetilde{H}_z(n),0\le z\le x\})=x$ a.s., and the distribution of $\widetilde{H}_x(n)$ is such that 
  $$\bE\left(e^{-\beta\widetilde{H}_x(n)}\right)=\exp\left\{-x\beta-x\int_{[0,2^{-n}]}(1-e^{-\beta y})\Lambda_H(dy)\right\}\rightarrow e^{-x\beta}\qquad\mbox{for all $\beta\ge0$.}$$
  Therefore, $\widetilde{H}_x(n)\rightarrow x$ in probability, and there is a subsequence $(n_k)_{k\ge 0}$ along which 
  convergence holds almost surely. Now let $\varepsilon>0$. Then there is (random) $N\ge 0$ such that for all $n_k\ge N$, we have
  $x\le \widetilde{H}_x(n_k)\le x+\varepsilon$. Therefore, 
  $$\limsup_{k\rightarrow\infty}\widetilde{K}_x(n_k)\le\inf_{\varepsilon>0}K_{x+\varepsilon}=K_x.$$ 
  Since $L_n(x)\le x$ a.s., the previous lemma implies the claimed convergence.
\end{pf} 

\begin{pfofpropP1} This proof is for the case where $X$ has a full right inverse and we only prove (\ref{P1E1}). The general case can be adapted. By
  Lemma \ref{lm9}, we can approximate $K_x=\lim_{k\rightarrow\infty}\widetilde{K}_x(n_k)$. The Laplace exponent of 
  $\widetilde{K}_x(n_k)$ follows from Lemma \ref{lm7} and this yields
  \beq\rho(q)&=&-\ln\left(\bE\left(e^{-qK_1}\right)\right)=-\lim_{k\rightarrow\infty}\ln\left(\bE\left(e^{-q\widetilde{K}_1(n_k)}\right)\right)\\
&=&\eta q+\lim_{k\rightarrow\infty}\int_{[0,\infty)\times[0,2^{-{n_k}}]}(1-e^{-qs})\Lambda(ds,dy)\\
&&\hspace{1cm}+\lim_{k\rightarrow\infty}\int_{[0,\infty)\times(2^{-{n_k}},\infty)}\int_{[0,\infty)\times[0,\infty)}(1-e^{-qt})\bP(s+T_{\{-y\}}\in dt,0\in dz)\Lambda(ds,dy)\\
&=&\eta q+\int_{[0,\infty)\times\{0\}}(1-e^{-qs})\Lambda(ds,dy)+\int_{[0,\infty)\times(0,\infty)}\left(1-e^{-qs}\bE\left(e^{-qT_{\{-y\}}}\right)\right)\Lambda(ds,dy).
  \eeq
\end{pfofpropP1}

\end{appendix}


\begin{thebibliography}{99}


\bibitem{bert96} Bertoin, J. (1996) \textit{L\'evy Processes}. Cambridge
University Press.




\bibitem{Doney-07} Doney, R (2007) {\it Fluctuation theory for L\'evy processes. Lectures from the 35th Summer School on Probability Theory held in Saint-Flour}. Lecture Notes in Mathematics \bf 1897\rm. Springer, Berlin.

\bibitem{DoK-06} Doney, R.A. and Kyprianou, A.E. (2006) Overshoots and undershoots of L\'evy processes. \it Ann. Appl. Prob. \rm\bf 16\rm, No. 1, 91--106. 

\bibitem{DoSa-10} Doney, R. and Savov, M. (2010) Right inverses of L\'evy processes. \it to appear in Ann. Prob. \rm


\bibitem{Ev-00} Evans, S. (2000) Right inverses of non-symmetric L\'evy processes and stationary stopped local times. \it Probab. Theory Related Fields. \rm\bf 118\rm, 37--48 


\bibitem{Sat-99} Sato, K. (1999) {\it L\'evy processes and infinitely divisible distributions}. Cambridge University Press.

\bibitem{Sim-99} Simon, T. (1999) Subordination in the wide sense for L\'evy processes. \it Probab. Theory Related Fields. \rm\bf 115\rm, 445--477

\bibitem{Win-02} Winkel, M. (2002) Right inverses of non-symmetric L\'evy processes. \it Ann. Prob. \rm\bf 30\rm, No. 1, 382--415.

\end{thebibliography}
\end{document}